\documentclass[english]{amsart}
\usepackage{mathptmx}
\usepackage{helvet}
\usepackage{beramono}
\usepackage[T1]{fontenc}
\usepackage[latin9]{inputenc}
\usepackage{listings}
\usepackage[b5paper]{geometry}
\usepackage{babel}
\usepackage{float}
\usepackage{amsthm}
\usepackage{amssymb}
\usepackage{graphicx}
\usepackage{esint}
\usepackage[unicode=true,pdfusetitle,
 bookmarks=true,bookmarksnumbered=false,bookmarksopen=false,
 breaklinks=false,pdfborder={0 0 1},backref=false,colorlinks=false]
 {hyperref}

\makeatletter

\providecommand{\tabularnewline}{\\}
\floatstyle{ruled}
\newfloat{algorithm}{tbp}{loa}
\providecommand{\algorithmname}{Algorithm}
\floatname{algorithm}{\protect\algorithmname}

\numberwithin{equation}{section}
\theoremstyle{plain}
\newtheorem{thm}{\protect\theoremname}
  \theoremstyle{definition}
  \newtheorem{defn}[thm]{\protect\definitionname}
  \theoremstyle{plain}
  \newtheorem{fact}[thm]{\protect\factname}
  \theoremstyle{definition}
  \newtheorem{example}[thm]{\protect\examplename}
  \theoremstyle{remark}
  \newtheorem{claim}[thm]{\protect\claimname}

\makeatletter
\renewcommand*{\@biblabel}[1]{\hfill#1.}
\makeatother

\makeatother

  \providecommand{\claimname}{Claim}
  \providecommand{\definitionname}{Definition}
  \providecommand{\examplename}{Example}
  \providecommand{\factname}{Fact}
\providecommand{\theoremname}{Theorem}

\begin{document}

\title{Partial freeness of random matrices}

\author{Jiahao Chen}

\address{Department of Chemistry, Massachusetts Institute of Technology, 77
Massachusetts Avenue, Cambridge MA 02139}

\email{\texttt{jiahao@mit.edu}}

\author{Troy Van Voorhis}

\address{Department of Chemistry, Massachusetts Institute of Technology, 77
Massachusetts Avenue, Cambridge MA 02139}

\email{\texttt{tvan@mit.edu}}

\author{Alan Edelman}

\address{Department of Mathematics, Massachusetts Institute of Technology,
77 Massachusetts Avenue, Cambridge MA 02139}

\email{\texttt{edelman@math.mit.edu}}
\begin{abstract}
We investigate the implications of free probability for finite-dimensional,
Hermitian random matrices. While most pairs of matrices are not free,
they can nevertheless be characterized by the joint moments which
violate freeness. We use this to extend the notion of freeness to
an intermediate property which we call partial freeness. Furthermore,
we can calculate the deviation of the true density of states from
the free convolution using asymptotic moment expansions, whose coefficients
have a natural combinatorial interpretation over weights of closed
paths. We have developed computer programs for characterizing partial
freeness from either numerical samples of matrices or analytic eigenvalue
distributions. Our results provide a rich extension of free probability
into the realm of finite linear algebra, with freeness emerging in
limiting cases. 
\end{abstract}
\maketitle
\tableofcontents{}

\section*{Introduction}

Free probability has received much attention since its discovery as
an algebraic structure for noncommuting op\-e\-ra\-tors.\cite{Voiculescu1985}
Subsequently, it has found a place in combinatorics with the deep
relationship between free cumulants and noncrossing partitions.\cite{Nica2006a,Novak2011}
Free probability for random matrices usually focuses on the asymptotic
freeness of infinite matrices.\cite{Voiculescu1991,Biane1998} In
contrast, we investigate here how free probability offers us new perspectives
on finite-dimensional matrices. We develop and extend the notion of
freeness to finite random matrices using linear algebra and elementary
statistics, without requiring intricate knowledge of operator algebras
or combinatorics.

In this paper, we consider the problem of calculating eigenvalues
of sums of finite matrices given the eigenvalues of the individual
matrices, as an illustration of the power of free probability theory.
In general, the eigenvalues of the sum of two matrices $A+B$ are
not simply the sums of the eigenvalues of the individual matrices
$A$ and $B$;\cite{Knutson2001} as matrices do not generally commute,
the addition of eigenvalues must take into account the relative orientations
of eigenvectors. However, free probability does allows us to do this
calculation in the limiting case where the rotation matrix between
the two bases is so random as to be uniformly oriented, i.e. of uniform
Haar measure. The matrices $A$ and $B$ are then said to be in generic
position, or free, and the eigenvalue spectrum of $A+B$ converges,
in a sense, to the additive free convolution $A\boxplus B$ of two
random matrices $A$ and $B$ as the matrix dimensions increase to
infinity.\cite{Nica2006a}

A natural question to ask is how accurately $A\boxplus B$ approximates
the exact eigenvalue spectrum, or density of states (d.o.s.), of the
sum $A+B$ when the individual matrices are known to be noncommuting
but not necessarily free. We seek to quantify this statement in this
paper. In addition to classifying two random matrices as being free
or not free relative to each other, we can characterize them as having
an intermediate, graduated property which we call partial freeness,
and furthermore we are able to quantify the leading-order discrepancy
between freeness and partial freeness. This has already helped us
explain the unexpected accuracy of approximations to the Hamiltonians
of disordered condensed matter systems.\cite{Chen2012}

We begin with a brief, self-contained review of free probability from
a random matrix theoretic perspective, and provide an elementary illustration
of how computing the additive free convolution using an integral transform
allows us to calculate the d.o.s.\ for the sum of free random matrices.
Next, we recap how the additive free convolution can also be approached
via the moments of random matrices, and in particular how both classical
and free independence can be interpreted as imposing precise rules
for the decomposition of joint moments of arbitrary orders. We then
show how we can generalize this to the notion of partial freeness
and describe a procedure for detecting it numerically from samples
of random matrices.

\section{Freeness of two matrices}

We use the notation $\left\langle \cdot\right\rangle $ for the normalized
expected trace (n.e.t.) $\frac{1}{N}\mathbb{E}\mbox{ Tr}\cdot$ of
a $N\times N$ matrix.
\begin{defn}
The random matrices $A$ and $B$ are free (or synonymously, freely
independent) with respect to the n.e.t.\ if for all $k\in\mathbb{N}$,
\begin{equation}
\left\langle p_{1}\!\left(A\right)q_{1}\!\left(B\right)p_{2}\!\left(A\right)q_{2}\!\left(B\right)\cdots p_{k}\!\left(A\right)q_{k}\!\left(B\right)\right\rangle =0\label{eq:free-defn}
\end{equation}
for all polynomials $p_{1},q_{1},p_{2},q_{2},\dots p_{k},q_{k}$ such
that $\left\langle p_{1}\!\left(A\right)\right\rangle =\left\langle q_{1}\!\left(B\right)\right\rangle =\cdots=0$.\cite[Definition 4.2]{Voiculescu1985}

This generalizes the notion of (classical) independence of scalar
random variables: were $A$ and $B$ to commute, the preceding with
$k=1$ would suffice.\end{defn}
\begin{fact}
The preceding is equivalent to defining free independence using the
special case of the centering polynomials $p_{i}\left(x\right)=x^{n_{i}}-\left\langle x^{n_{i}}\right\rangle $,
$q_{i}\left(x\right)=x^{m_{i}}-\left\langle x^{m_{i}}\right\rangle $,
$i=1,\dots,k$ for positive integers $n_{1},m_{1},\dots,n_{k},m_{k}$.\cite[Proposition 4.3]{Voiculescu1985}
\end{fact}
That this is sufficient follows from the linearity of the n.e.t. In
principle, this would allow us to check if two matrices $A$ and $B$
were free by checking that all centered joint moments 
\[
\left\langle \left(A^{n_{1}}\!-\!\left\langle A^{n_{1}}\right\rangle \right)\left(B^{m_{1}}\!-\!\left\langle B^{m_{1}}\right\rangle \right)\cdots\left(A^{n_{k}}\!-\!\left\langle A^{n_{k}}\right\rangle \right)\left(B^{m_{k}}\!-\!\left\langle B^{m_{k}}\right\rangle \right)\right\rangle 
\]
vanish for all positive exponents $n_{1}$,$m_{1}$,$\ldots$,$n_{k}$,$m_{k}$.
However, this is numerically impractical due to the need to check
joint moments of all orders, as well as the presence of fluctuations
from sampling error if using Monte Carlo, which causes the higher
order joint moments to converge slowly. In practice, it is far easier
to check for freeness by examining how the d.o.s.\ of the exact sum
$A+B$ converges to the p.d.f.\ defined by the free convolution $A\boxplus B$,\cite{Voiculescu1986}
which we will now define.

\subsection{The free convolution}
\begin{defn}
The $R$-transform of the p.d.f. $f_{A}$, denoted by $R_{A}$, is
defined implicitly via the following Cauchy trans\-form:\cite{Voiculescu1985,Nica2006a}%
\footnote{There are unfortunately two extant notations for the $R$-transform.
We use here the $R$-transform as presented in \cite{Nica2006a};
this differs slightly from the original notation of Voiculescu,\cite{Voiculescu1985}
which is the $\mathcal{R}$-transform elsewhere, e.g.\ in \cite{Nica2006a}.
The relationship between the two is $R\left(w\right)=w\mathcal{R}\left(w\right)$.%
}
\begin{equation}
w=\lim_{\epsilon\downarrow0}\int_{\mathbb{R}}\frac{f_{A}\left(z\right)}{R_{A}\left(w\right)-\left(z+i\epsilon\right)}dz.\label{eq:r-transform}
\end{equation}
Some intuition for the $R$-transform may be achieved by expanding
the Cauchy integral as a formal power series:%
\footnote{Physicists may recognize $G_{A}\left(w\right)$ as the retarded Green
function corresponding to the Hamiltonian $A$.%
}
\begin{equation}
G_{A}\left(w\right)=\lim_{\epsilon\downarrow0}\int_{\mathbb{R}}\frac{f_{A}\left(z\right)}{w-z-i\epsilon}dz=\sum_{k=0}^{\infty}\frac{\mu_{k}\left(A\right)}{w^{k+1}},
\end{equation}
where $\mu_{k}$ is the $k$th moment
\begin{equation}
\mu_{k}\left(A\right)=\int_{\mathbb{R}}x^{k}f_{A}\left(x\right)dx=\left\langle A^{k}\right\rangle .
\end{equation}
In other words, the Cauchy transform of a p.d.f.\ is a generating
function of its moments. We then have that the $R$-transform $R_{A}$
inverts the Cauchy transform $G_{A}$ in the functional sense, i.e.\ that
\begin{equation}
G_{A}\left(R_{A}\left(w\right)\right)=w.
\end{equation}
Viewing both $G_{A}$ and $R_{A}$ as formal power series, the latter
is simply the reversion of the former,\cite{Morse1953,Henrici1974}
in the sense that $R_{A}$ is a series in $w$ whose inverse with
respect to composition is $G_{A}$ as a series in $1/z$. The coefficients
of the $R$-transform are then the free cumulants $\nu_{k}$, i.e.
\end{defn}
\begin{equation}
R_{A}\left(w\right)=\sum_{k=0}^{\infty}\nu_{k}w^{k-1},
\end{equation}
with $\nu_{0}=1$. The free cumulants are particular combinations
of moments $\nu_{k}=\nu_{k}\left(\mu_{1},\dots,\mu_{k}\right)$ which
shall be made more explicit later.
\begin{defn}
The free convolution $A\boxplus B$ is defined via its $R$-transform
\begin{equation}
R_{A\boxplus B}\left(w\right)=R_{A}\left(w\right)+R_{B}\left(w\right)-\frac{1}{w}.
\end{equation}

\end{defn}
The free cumulants linearize the free convolution in the sense that
for all $k>0$,\cite{Nica2006a,Novak2011}
\begin{equation}
\nu_{k}\left(A\boxplus B\right)=\nu_{k}\left(A\right)+\nu_{k}\left(B\right),
\end{equation}
and the subtraction of $1/w$ produces a properly normalized p.d.f.
by conserving $\nu_{0}\left(A\boxplus B\right)=\mu_{0}\left(A\boxplus B\right)=1$.

In Section~\ref{sub:ex-free-finite}, we show an example of calculating
$f_{A\boxplus B}$ analytically via the $R$-transform. In general,
such analytic calculations are hindered by the functional inversions
required in (\ref{eq:r-transform}). This has inspired interesting
work in calculating $A\boxplus B$ numerically, such as in the RMTool
package.\cite{Rao2007,Olver2012} We discuss instead an alternate
strategy starting directly from numerical samples of random matrices,
which generalizes naturally to general pairs of matrices. In situations
where only the numerical samples are known, it may be convenient instead
to use the result of Fact~\ref{A-QBQ} described in the next section.

\subsection{Free convolution from random rotations}
\begin{defn}
A square matrix $Q$ is a unitary/\-orthogonal/\-symplectic random
matrix of Haar measure if for any constant unitary/orthogonal/symplectic
matrix $P$, the integral of any function over $dQ$ is identical
to the integral over $d\left(PQ\right)$ or that over $d\left(QP\right)$.\end{defn}
\begin{example}
Unitary matrices of dimension $N=1$ are simply scalar unit complex
phases of the form $e^{i\theta}$. Haar measure over $e^{i\theta}$
can be written simply as $\mbox{d}\theta/2\pi$. This is manifestly
rotation invariant, as multiplying $e^{i\theta}$ by any constant
phase factor $e^{i\phi}$ simply changes the measure to $\mbox{d}(\theta+\phi)/2\pi=\mbox{d}\theta/2\pi$.
\end{example}
Uniform Haar measure generalizes the concept of uniformity to higher
dimensions by preserving the notion of invariance with respect to
arbitrary rotations. Consequently, the eigenvalues of $Q$ lie uniformly
on the unit circle on the complex plane.\cite{Diaconis1994} Explicit
samples can be generated numerically by performing $QR$ decompositions
on $N\times N$ matrices sampled from the Gaussian orthogonal (unitary)
ensemble.\cite{Diaconis2005}
\begin{fact}
\label{A-QBQ}For a pair of Hermitian (real symmetric) random matrices
$A$ and $B$, the d.o.s.\ of $A+QBQ^{\dagger}$, where $Q$ is a
unitary (orthogonal) random matrix of Haar measure, coincides with
the p.d.f.\ of $A\boxplus B$ in the limit of infinitely large matrices
$N\rightarrow\infty$.
\end{fact}
Consider the diagonalization of $A=Q_{A}\Lambda_{A}Q_{A}^{\dagger}$
and $B=Q_{B}\Lambda_{B}Q_{B}^{\dagger}$. The d.o.s.\ of $A+QBQ^{\dagger}$
is identical to that of $\Lambda_{A}+\left(Q_{A}^{\dagger}QQ_{B}\right)\Lambda_{B}\left(Q_{B}^{\dagger}Q^{\dagger}Q_{A}\right)$,
since these matrices are related by the similarity transformation
$Q_{A}^{\dagger}\left(\cdot\right)Q_{A}$. However, the Haar property
of $Q$ means that the d.o.s.\ of this matrix is identical to that
of $\Lambda_{A}+Q\Lambda_{B}Q^{\dagger}$. This gives us another interpretation
of free convolution: it describes the statistics resulting from adding
two random matrices when the basis of one matrix is randomly rotated
or ``spun around'' relative to the other. The information about
the relative orientations of the two bases is effectively ignored,
retaining only the knowledge that they are not parallel so that $A$
and $B$ do not commute.

The freeness of random matrices usually discussed only in the limit
of infinitely large matrices, where it is called asymptotic freeness.\cite{Voiculescu1991}
For example, two matrices sampled from the Gaussian ensembles (orthogonal,
unitary or symplectic) are free.\cite{Nica2006a} Nevertheless, finite-dimensional
random matrices can exhibit freeness as well. We now provide some
examples and illustrate the analytic calculation of the free convolution
using the $R$-transform.

\subsection{Examples of free finite-dimensional matrices\label{sub:ex-free-finite}}
\begin{figure}[h]
\caption{\label{fig:arcsine-vs-binomial}The d.o.s. $f_{A\boxplus B}\left(x\right)$
and $f_{A\star B}$$\left(x\right)$ for the free (dashed blue line)
and classical convolutions (solid black bars) of the matrices in Example
1, as given in (\ref{eq:arcsine}) and (\ref{eq:discrete-binomial})
respectively. The heights of the lines in the plot of $f_{A\star B}$
indicate the point masses.}

\includegraphics[width=0.95\columnwidth,height=0.95\columnwidth,keepaspectratio]{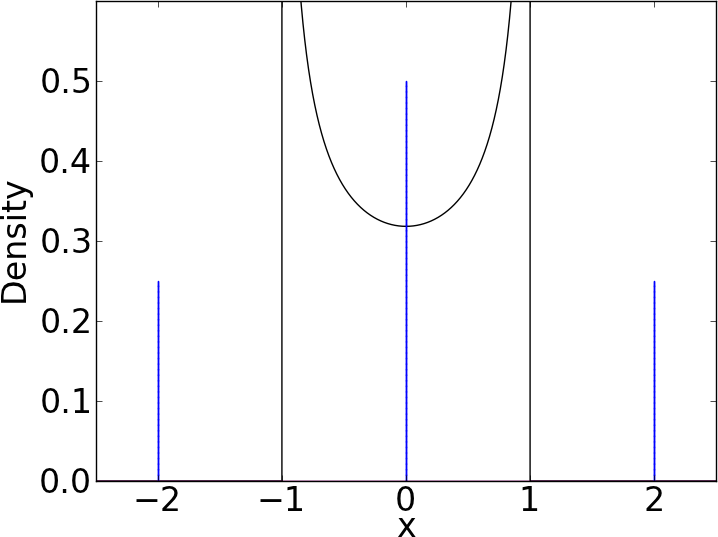}
\end{figure}

\begin{example}
Consider the $2\times2$ real symmetric random matrices 
\begin{equation}
A\left(t\right)=U\left(t\right)\sigma_{z}U\left(-t\right),\quad B\left(t\right)=U\left(-t\right)\sigma_{z}U\left(t\right),
\end{equation}
where $\sigma_{z}$ is the Pauli matrix $\left(\begin{array}{cc}
1 & 0\\
0 & -1
\end{array}\right)$, $U\left(t\right)$ is the rotation matrix $\left(\begin{array}{cc}
\cos t & \sin t\\
-\sin t & \cos t
\end{array}\right)$, and the rotation angle $t$ is uniformly sampled on the interval
$\left[0,\pi\right)$. By construction, the d.o.s.\ of $A\left(t\right)$
and $B\left(t\right)$ are identical; their eigenvalues have the p.d.f.
\begin{equation}
f_{A}\left(x\right)=f_{B}\left(x\right)=\frac{1}{2}\left(\delta\left(x+1\right)+\delta\left(x-1\right)\right),
\end{equation}
where $\delta\left(x\right)$ is the Dirac delta distribution. Furthermore,
for any particular $t$, the sum of $A\left(t\right)$ and $B\left(t\right)$
can be written in the basis where $A\left(t\right)$ is diagonal as
\begin{equation}
M\left(t\right)=\sigma_{z}+U\left(-2t\right)\sigma_{z}U\left(2t\right).
\end{equation}
By construction, $U\left(2t\right)$ is of uniform Haar measure and
so the d.o.s.\ of $M\left(t\right)$ is given exactly by the additive
free convolution of $A\left(t\right)$ and $B\left(t\right)$.\cite{Voiculescu1991,Nica2006a}
The $R$-transforms of $f_{A}$ and $f_{B}$ are 
\begin{equation}
R_{A}\left(w\right)=R_{B}\left(w\right)=\frac{1\pm\sqrt{1+4w^{2}}}{2w}.
\end{equation}
Performing the free convolution,
\begin{equation}
R_{A\boxplus B}\left(w\right)=R_{A}\left(w\right)+R_{B}\left(w\right)-\frac{1}{w}=\pm\frac{\sqrt{1+4w^{2}}}{w}.
\end{equation}
Finally, we calculate the p.d.f.\ using the Plemelj inversion formula:\begin{subequations}
\begin{align}
f_{A\boxplus B}\left(x\right) & =\frac{1}{\pi}\left[\mbox{Im }R_{A\boxplus B}^{-1}\left(w\right)\right]_{w=x}\label{eq:plemelj}\\
 & =\frac{1}{\pi\sqrt{4-x^{2}}},\label{eq:arcsine}
\end{align}
\end{subequations}which is the arcsine distribution on the interval
$\left[-2,2\right]$, and we have retained only the positive root
to obtain a nonnegative probability density. The odd moments all vanish
by the even symmetry of $f_{\negthinspace A\boxplus B}$, and the
even moments are the central binomial coefficients $\mu_{2n}\negthinspace\left(\negthinspace A\negthinspace\boxplus\negthinspace B\right)$
$=\binom{2n}{n}$.
\end{example}
This example shows that the free convolution of two discrete probability
distributions can be a continuous probability distribution. In contrast,
the classical convolution $A\star B$ produces the p.d.f.
\begin{equation}
f_{A\star B}\left(x\right)=f_{A}\star f_{B}=\int_{\mathbb{R}}f_{A}\left(y\right)f_{B}\left(x-y\right)dy=\frac{1}{4}\left(\delta\left(x+2\right)+2\delta\left(x\right)+\delta\left(x-2\right)\right),\label{eq:discrete-binomial}
\end{equation}
which is simply a discrete binomial distribution. The results of the
two convolutions are plotted in Figure~\ref{fig:arcsine-vs-binomial}.

Here is another example of matrices with the same d.o.s.\ as before.
This time, the matrices are not random at all.
\begin{example}
\label{example:pauli-sx}The $2N\times2N$ deterministic matrices
\begin{equation}
A=\left(\begin{array}{ccccc}
0 & 1\\
1 & 0\\
 &  & 0 & 1\\
 &  & 1 & 0\\
 &  &  &  & \ddots
\end{array}\right),\quad B=\left(\begin{array}{ccccc}
0 &  &  &  & 1\\
 & 0 & 1\\
 & 1 & 0\\
 &  &  & \ddots\\
1 &  &  &  & 0
\end{array}\right)
\end{equation}
are asymptotically free as $N\rightarrow\infty$. Each consists of
$N$ direct sums of the Pauli matrix $\sigma_{x}=\left(\begin{array}{cc}
0 & 1\\
1 & 0
\end{array}\right)$, with $B$ having a basis shifted relative to $A$ with circulant
(periodic) boundary conditions. The d.o.s.\ are the same as in the
previous example and so the calculations of $A\boxplus B$ and $A\star B$
proceed identically. Considering the matrix $B^{\prime}$, being $B$
without the circulant entries on the lower left and upper right, we
also have that $A$ and $B^{\prime}$ are asymptotically free as $N\rightarrow\infty$.
\end{example}

\subsection{Comparison of free and classical convolutions}

To conclude this introductory survey, we compare the free additive
convolution $\boxplus$ and the classical convolution $\star$ in
more detail. We note that the Fourier transform $\;\widehat{\cdot}\;$
turns classical convolutions into products according to the convolution
theorem, i.e. $\widehat{f_{A}\star f_{B}}=\widehat{f_{A}}\widehat{f_{B}}$;
taking the logarithm allows this to be written as a linear sum:
\begin{equation}
\log\widehat{f_{A}\star f_{B}}=\log\widehat{f_{A}}+\log\widehat{f_{B}}.
\end{equation}

Furthermore if $f_{A}$ and $f_{B}$ are p.d.f.s, then $\log\widehat{f_{A}}$
and $\log\widehat{f_{B}}$ can be identified as the corresponding
(classical) cumulant generating functions, i.e. $\log\widehat{f_{A}}$
can be expanded in a formal power series
\begin{equation}
\log\widehat{f_{A}}\left(w\right)=\sum_{k=0}^{\infty}\frac{\kappa_{n}\left(A\right)w^{n}}{n!},
\end{equation}
where $\kappa_{n}\left(A\right)$ is the $n$th (classical) cumulant
of $f_{A}$, and similarly for $f_{B}$. We have also have that $\kappa_{n}\left(A\star B\right)=\kappa_{n}\left(A\right)+\kappa_{n}\left(B\right)$
for $n\ge1$, i.e. that cumulants linearize the convolution. Drawing
an analogy between the cumulants $\kappa_{n}$ and the free cumulants
$\nu_{n}$, the $R$-transform is often described as the free analogue
of the log-Fourier transform, with the Cauchy transform being the
analogue of the Fourier transform and the functional inversion playing
the part analogous to the logarithm.

Finally, we note that the sum of two scalar random variables $a$
and $b$, sampled with p.d.f.s $f_{A}$ and $f_{B}$ respectively,
itself has the resulting p.d.f.\ $f_{A\star B}$.\cite{Feller1971}
For matrices $A$ and $B$ with d.o.s. $f_{A}$ and $f_{B}$ respectively,
then the matrix $A+\Pi B\Pi^{T}$, where $\Pi$ is a random permutation
matrix, has d.o.s.\ $f_{A}\star f_{B}$. This is equivalent to the
p.d.f.\ formed by picking an eigenvalue of $A$ at random and an
eigenvalue of $B$ at random. In this sense, the discrete random permutation
$\Pi$ which generates the classical convolution is the analogue of
the continuous random rotation $Q$ of uniform Haar measure in free
convolution.

Table~\ref{tab:free-vs-classical} summarizes the analogies between
the free and classical convolutions.

\begin{table}[h]
\caption{\label{tab:free-vs-classical}Correspondence between the free and
classical convolutions.}

\begin{tabular}{|c|c|}
\hline 
$A\boxplus B$ & $A\star B$\tabularnewline
\hline 
\hline 
$R$-transform $R$ & log-Fourier transform $\log\hat{f}$\tabularnewline
\hline 
Cauchy transform $G$ & Fourier transform $\hat{f}$\tabularnewline
\hline 
functional inversion & logarithm\tabularnewline
\hline 
free cumulants $\nu_{n}$ & (classical) cumulants $\kappa_{n}$\tabularnewline
\hline 
Plemelj inversion & inverse Fourier transform\tabularnewline
\hline 
\hline 
$A+QBQ^{\dagger}$ & $A+\Pi B\Pi^{T}$\tabularnewline
\hline 
uniform Haar measure $Q$ & random permutations $\Pi$\tabularnewline
\hline 
\end{tabular}
\end{table}

\section{Density of states of sums of random matrices}

The introductory examples show that if $A$ and $B$ are free, then
the d.o.s.\ of $A+B$ can be calculated exactly without any detailed
knowledge of the eigenvectors of $A$ and $B$. However, not all pairs
of random matrices are free. Nevertheless, the free convolution can
provide surprisingly accurate approximations to their d.o.s.\ even
in the general case, and the degree of approximation can even be quantified
by examining individual moments of the sum, $\mu_{n}\left(A+B\right)$
as shown in Section~\ref{partial-freeness}. In order to do this,
however, we will first need to examine how each moment of $A+B$ subdivides
into sums over joint moments of $A$ and \textbf{$B$}, and how the
primordial definition of freeness in (\ref{eq:free-defn}) determines
what these moments must be for free $A$ and $B$. For simplicity,
we assume in this paper that all necessary moments $\left\{ \mu_{n}\right\} $
exist, so that the moments capture all the information contained in
the corresponding p.d.f. For random matrices, this simply means that
all powers of the matrix must have a finite n.e.t.

\subsection{Calculating moments of $A+B$ from the joint moments of $A$ and
$B$}

For general $A$ and $B$, it is possible to characterize the d.o.s.\ of
$A+B$ completely without constructing the sum and diagonalizing it
if all their joint moments are known. We can then calculate all the
moments $\left\{ \mu_{n}\left(A+B\right)\right\} $ from the definition
of the moment of a random matrix:\begin{subequations} 
\begin{align}
\mu_{n}= & \left\langle \left(A+B\right)^{n}\right\rangle \label{eq:moment-def}\\
= & \langle A^{n}+A^{n-1}B+A^{n-2}BA+\cdots+BA^{n-1}+A^{n-2}B^{2}+\cdots+B^{2}A^{n-2}+\cdots+B^{n}\rangle\label{eq:moment-binomial}\\
= & \left\langle A^{n}\right\rangle +n\left\langle A^{n-1}B\right\rangle +n\left\langle A^{n-2}B^{2}\right\rangle +\cdots+\left\langle B^{n}\right\rangle .\label{eq:moment-expanded}
\end{align}
\end{subequations}The second equality follows directly from the noncommutative
binomial expansion of (\ref{eq:moment-binomial}), and the third equality
follows from the linearity of $\left\langle \cdot\right\rangle $
and its cyclic invariance, i.e.\ $\left\langle AB\right\rangle =\left\langle BA\right\rangle $.

We refer to (\ref{eq:moment-expanded}) as the \emph{word expansion
of $\mu_{n}$}. As written, there are $2^{n}$ terms in (\ref{eq:moment-binomial})
but some of them yield the same n.e.t.\ in (\ref{eq:moment-expanded})
identically because of cyclic invariance. The equivalence classes
defined by grouping identical terms in this manner are exactly those
of combinatorial necklaces.\cite{MacMahon1892,Riordan1957}
\begin{defn}
\label{def:necklace}An $\left(n,k\right)$-word $W$ is a string
of $n$ symbols, each of which can have any of $k$ values. An $\left(n,k\right)$-necklace
$\left[\mathcal{N}\right]$ is the equivalence class over $\left(n,k\right)$-words
$W$ with respect to cyclic permutations $\Pi$ of length $n$, i.e.
\begin{equation}
\left[\mathcal{N}\right]=\left\{ w\in W\vert\exists\pi\in\Pi:\;\mathcal{N}=\pi w\right\} .
\end{equation}

\end{defn}
There are efficient algorithms for enumerating all $\left(n,k\right)$-necklaces
for a given $n$ and $k$.\cite{Ruskey1992,Sawada2001} Furthermore,
the total number of terms in the word expansion (\ref{eq:moment-expanded})
is well-known:
\begin{fact}
The number of $\left(n,k\right)$-necklaces is
\begin{equation}
N\left(n,k\right)=\frac{1}{n}\sum_{d\vert n}\phi\left(d\right)k^{n/d}=\frac{1}{n}\sum_{i=1}^{n}k^{\gcd\left(i,n\right)},
\end{equation}
where $d\vert n$ means that $d$ divides $n$, $\phi$ is the Euler
totient function, and $\gcd$ is the greatest common divisor.\cite{MacMahon1892,Riordan1957}
By definition, $\phi\left(d\right)$ is the number of integers $m$
in the range $1\le m\le d$ that are relatively prime to $d$, i.e.
$gcd\left(d,m\right)=1$.
\end{fact}
In addition, we can determine the multiplicity of each term in (\ref{eq:moment-expanded}),
which is identical to the number of words in the equivalence class
defined by each corresponding particular necklace. We state this very
simple fact without proof and provide an example.
\begin{fact}
Let $m=\#\left(\left[\mathcal{N}\right]\right)$ be the number of
$\left(n,k\right)$-words belonging to the equivalence class that
defines the necklace $\left[\mathcal{N}\right]$. Then $m$ is the
length of the longest cyclic permutation that leaves any word $W\in\mathcal{N}$
unchanged, i.e. it is the length of the longest subword $S$ of a
word $W\in\mathcal{N}$ such that $W=S^{n/m}$.\end{fact}
\begin{example}
The necklace $\left[AABAAB\right]=\left[A^{2}BA^{2}B\right]$ is an
equivalence class over $\left(6,2\right)$-words of size 3, since
applying a (one-symbol) cyclic permutation three times leaves $A^{2}BA^{2}B$
unchanged:
\begin{equation}
AABAAB\mapsto ABAABA\mapsto BAABAA\mapsto AABAAB,
\end{equation}
i.e. $\#\left(\left[A^{2}BA^{2}B\right]\right)=3$ which follows from
the fact that $AABAAB=\left(AAB\right)^{2}$.
\end{example}
The algorithmic enumeration of necklaces and their multiplicities
allow us to sum the joint moments in the word expansion (\ref{eq:moment-expanded})
to obtain $\mu_{n}$. As $N\left(n,k\right)=\mathcal{O}\left(k^{n}/n\right)$
asymptotically as $n\rightarrow\infty$, the word expansion saves
approximately a factor of $n$ in effort relative to working with
the naive noncommutative binomial expansion, which has $k^{n}$ terms.

\subsection{Decomposition rules for joint moments}

We have reduced the problem of calculating $\mu_{n}$ to that of calculating
joint moments; each has the form $\left\langle A^{n_{1}}B^{m_{1}}\cdots A^{n_{k}}B^{m_{k}}\right\rangle $
for positive integers $n_{1}$, $m_{1}$, $\ldots$, $n_{k}$, $m_{k}$.
In general, this is not the most compact way to specify the relationship
between $A$ and $B$. However, classical and free independence each
provide a prescription for computing such joint moments in terms of
the pure moments $\left\langle A\right\rangle $, $\left\langle A^{2}\right\rangle $,
$\ldots$ and $\left\langle B\right\rangle $, $\left\langle B^{2}\right\rangle $,
$\ldots$ of $A$ and $B$ respectively.
\begin{fact}
For classically independent random matrices $A$ and $B$, 
\begin{equation}
\left\langle A^{n_{1}}B^{m_{1}}\cdots A^{n_{k}}B^{m_{k}}\right\rangle =\left\langle A^{n_{1}+\cdots+n_{k}}B^{m_{1}+\cdots+m_{k}}\right\rangle =\left\langle A^{n_{1}+\cdots+n_{k}}\right\rangle \left\langle B^{m_{1}+\cdots+m_{k}}\right\rangle ,\label{eq:class-indep}
\end{equation}
i.e. $A$ and $B$ behave as if they commute.\cite{Nica2006a}
\end{fact}
The analogous rule for free independence is more complicated; however,
an implicit formula can be derived from the primordial definition
of freeness in (\ref{eq:free-defn}) by using the linearity of the
n.e.t.:\begin{subequations}
\begin{align}
0= & \left\langle \left(A^{n_{1}}-\left\langle A^{n_{1}}\right\rangle \right)\left(B^{m_{1}}-\left\langle B^{m_{1}}\right\rangle \right)\left(A^{n_{k}}-\left\langle A^{n_{k}}\right\rangle \right)\left(B^{m_{k}}-\left\langle B^{m_{k}}\right\rangle \right)\right.\\
= & \left\langle A^{n_{1}}B^{m_{1}}\cdots A^{n_{k}}B^{m_{k}}\right\rangle \nonumber \\
 & -\left\langle A^{n_{1}}\right\rangle \left\langle A^{n_{2}}B^{m_{2}}\cdots A^{n_{k}}B^{m_{k}+m_{1}}\right\rangle -\left\langle B^{m_{1}}\right\rangle \left\langle A^{n_{1}+n_{2}}B^{m_{2}}\cdots A^{n_{k}}B^{m_{k}}\right\rangle \nonumber \\
 & +\dots+\left(-1\right)^{\sum_{i=1}^{k}\left(n_{i}+m_{i}\right)}\left\langle A^{n_{1}}\right\rangle \cdots\left\langle A^{n_{k}}\right\rangle \left\langle B^{m_{1}}\right\rangle \cdots\left\langle B^{m_{k}}\right\rangle ,\label{eq:free-expandjoint}
\end{align}
\end{subequations}which can be rearranged immediately to give a recurrence
relation for the joint moment $\left\langle A^{n_{1}}B^{m_{1}}\cdots A^{n_{k}}B^{m_{k}}\right\rangle $
in terms of joint moments of lower order.

\section{\label{sec:Partial-freeness}Partial freeness}

The main result of our paper is to show that the following notion
of partial freeness is a useful generalization of freeness, particularly
for finite random matrices.
\begin{defn}
\label{partial-freeness} Two random matrices $A$ and $B$ are partially
free to order $p$ if the first difference between $A+B$ and $A\boxplus B$
occurs at the $p$th moment $\mu_{p}$, i.e. (\ref{eq:free-expandjoint})
holds for all joint moments of the form 
\[
\left\langle A^{n_{1}}B^{m_{1}}\cdots A^{n_{k}}B^{m_{k}}\right\rangle 
\]
with positive integers $n_{1}$, $m_{1}$, $\dots$, $n_{k}$, $\sum_{i=1}^{n}\left(n_{i}+m_{i}\right)=q$,
for all $q<p$, but there exists at least one joint moment for $q=p$
for which (\ref{eq:free-expandjoint}) does not hold. We say that
$A$ and $B$ are free to $p$ moments.
\end{defn}
In numerical applications, the difference between $A+B$ and $A\boxplus B$
must be tested for statistical significance if the joint moments are
calculated from Monte Carlo samples of $A$ and $B$.

This definition immediately allows us to restate the matching three
moments theorem of Refs.~\cite{Movassagh2010,Movassagh2011a}:
\begin{fact}
Let $A$ and $B$ be a pair of $N\times N$ diagonalizable random
matrices with $A=Q_{A}\Lambda_{A}Q_{A}^{\dagger}$ and $B=Q_{B}\Lambda_{B}Q_{B}^{\dagger}$.
If $\mathbb{E}\left[\left(Q_{B}^{\dagger}Q_{A}\right)_{ij}\right]=1/N$
for each matrix element of $Q_{B}^{\dagger}Q_{A}$, then $A$ and
$B$ are free to $p>3$ moments.
\end{fact}
Our definition is a natural generalization of the concept of freeness,
and they coincide if all the moments match.
\begin{claim}
Two random matrices $A$ and $B$ are free if they are partially free
to all orders.
\end{claim}
This follows immediately from the definitions of freeness and partial
freeness, so long as the limit $N\rightarrow\infty$ exists.

The following example illustrates that (partial) freeness can also
be a property of a pair of random and deterministic matrices.
\begin{example}
\label{example:tridiagonal}The $N\times N$ random matrix $A$, a
diagonal matrix with elements i.i.d.\ standard Gaussian random variates,
and $B$, the tridiagonal matrix
\[
\left(\begin{array}{cccc}
0 & 1 &  & 0\\
1 & \ddots & \ddots\\
 & \ddots & \ddots & 1\\
0 &  & 1 & 0
\end{array}\right)
\]
are partially free of order 8. This can be verified by explicit calculation
of the first eight moments. Again by even symmetry all the odd moments
of $A+B$ vanish while even moments are 1, 3, 17, 125, 1099, 11187,
129759,$\dots$ Furthermore we can identify the leading order deviation
as arising from the term $\left\langle \left(AB\right)^{4}\right\rangle =1$.
The significance of this result for condensed matter physics is discussed
in Ref.~\cite{Chen2012}.
\end{example}
In fact, partial freeness can also be a property of purely deterministic
matrices. Revisiting Example~\ref{example:pauli-sx}, we can show
that the $2N\times2N$-dimensional matrices $A$ and $B$ in that
Example are partially free to $2N$ moments. Since $A$ and $B$ are
each constructed out of direct sums of the same Pauli matrix, $A^{2}=B^{2}=I$
where $I$ is the identity matrix. Therefore, the demonstration of
partial freeness reduces to finding a $k$ for which $\mbox{tr }\left(AB\right)^{k}\ne0$.
As an illustration of what happens, consider that for $N=6$, we have
the sequence of matrices
\begin{equation}
\left\{ \left(AB\right)^{k}\right\} _{k=0}^{3}=\left\{ I,\left(\begin{array}{cccccc}
0 & 0 & 1 & 0 & 0 & 0\\
0 & 0 & 0 & 0 & 0 & 1\\
0 & 0 & 0 & 0 & 1 & 0\\
0 & 1 & 0 & 0 & 0 & 0\\
1 & 0 & 0 & 0 & 0 & 0\\
0 & 0 & 0 & 1 & 0 & 0
\end{array}\right),\left(\begin{array}{cccccc}
0 & 0 & 0 & 0 & 1 & 0\\
0 & 0 & 0 & 1 & 0 & 0\\
1 & 0 & 0 & 0 & 0 & 0\\
0 & 0 & 0 & 0 & 0 & 1\\
0 & 0 & 1 & 0 & 0 & 0\\
0 & 1 & 0 & 0 & 0 & 0
\end{array}\right),I\right\} .
\end{equation}
As $k$ increments by 1, the 1s in the odd columns move up two rows
(with wraparound) and the 1s in the even columns move down two rows.
Thus $k=N$ is the smallest positive integer for which the trace of
$\left(AB\right)^{k}$ does not vanish, and $A$ and $B$ are partially
free to $2N$ moments. We then recover the asymptotic freeness of
$A$ and $B$ immediately by taking the $N\rightarrow\infty$ limit.

Returning to the problem of computing the d.o.s.\ of the sum $A+B$,
we seek to ask how good an approximation the free convolution $A\boxplus B$
is to the d.o.s.\ of the sum $A+B$ when $A$ and $B$ are not free,
but only partially free. We now quantify this statement using asymptotic
moment expansions.\cite{Chen2012}

\section{Distinguishing between two distributions using asymptotic moment
expansions}

We have described partial freeness in terms of how the moments of
the sum $A+B$ differ from what free probability requires it to be.
This suggests that asymptotic moment expansions,\cite{Wallace1958}
which expand a p.d.f.\ $f$ about a reference p.d.f.\ $\tilde{f}$
and are parameterized by the moments (or cumulants) of the two distributions
being compared, provide a natural framework for examining how the
exact p.d.f.\ differs from the free convolution. We develop this
notion using the two standard expansions, namely the Gram--Charlier
series (of Type A) and the Edgeworth series.\cite{Stuart1994}

\subsection{The Gram--Charlier series}

The Gram--Charlier series arises immediately from the orthogonal polynomial
expansion with respect to $\tilde{f}$ as the weight:\cite[Chapter IX]{Szego1975}
\begin{equation}
f\left(x\right)=\sum_{n=0}^{\infty}c_{n}\phi_{n}\left(x\right)\tilde{f}\left(x\right),
\end{equation}
where the coefficients can be shown, by the orthonormality of the
orthogonal polynomials, to be
\begin{align}
\int_{\mathbb{R}}\phi_{m}\left(x\right)f\left(x\right)dx & =\sum_{n=0}^{\infty}c_{n}\int_{\mathbb{R}}\phi_{m}\left(x\right)\phi_{n}\left(x\right)\tilde{f}\left(x\right)dx=c_{m},
\end{align}
i.e.\ the $m$th coefficient is the expected value of the $m$th
orthogonal polynomial with respect to the probability density $f$.
By expressing the orthogonal polynomials in the monomial basis,
\begin{equation}
\phi_{m}\left(x\right)=\sum_{k=0}^{m}a_{mk}x^{k},\quad c_{m}=\sum_{k=0}^{m}a_{mk}\int_{\mathbb{R}}x^{k}f\left(x\right)dx=\sum_{k=0}^{m}a_{mk}\mu_{k},
\end{equation}
we get an explicit expansion of the Gram--Charlier coefficients $\left\{ c_{m}\right\} $
as linear combinations of the moments $\left\{ \mu_{k}\right\} $
of $f$. The so--called Gram--Charlier Type A series%
\footnote{This is often referred to simply as the Gram--Charlier series; however,
in this paper we mean the latter to be the generalization of the commonly
used Type A series to possibly non-Gaussian weight functions $f$.%
} is simply the special case of a standard Gaussian weight:
\begin{equation}
\tilde{f}\left(x\right)=\Phi\left(x\right)=\frac{1}{\sqrt{2\pi}}\exp\left(-\frac{x^{2}}{2}\right).
\end{equation}
The corresponding orthogonal polynomials $\left\{ \phi_{n}\right\} _{n}$
are the (pr\-obabilist's) Hermite polynomials $\left\{ He_{n}\right\} _{n}$.

\subsection{Edgeworth series}

The Gram--Charlier series can be seen as the output of an operator
$T$: 
\begin{equation}
T:\tilde{f}\rightarrow f,\quad T\left(x\right)=\sum_{n=0}^{\infty}c_{n}\phi_{n}\left(x\right),
\end{equation}
as applied to the reference p.d.f.\ $\tilde{f}$. In contrast, the
Edgeworth series is derived by rewriting $T$ as a differential operator,
as can be derived using the relations between a probability density,
its characteristic function $\chi\left(x\right)$ and the moment generating
function, and its cumulant generating function:
\begin{equation}
\chi\left(t\right)=\int_{\mathbb{R}}e^{itx}f\left(x\right)dx=\mathbb{E}_{f\left(x\right)}\left(e^{itx}\right)=\sum_{n=0}^{\infty}\frac{\mu_{n}}{n!}\left(it\right)^{n}=\exp\left(\sum_{n=1}^{\infty}\frac{\kappa_{n}}{n!}\left(it\right)^{n}\right).
\end{equation}
Writing down the analogous relations for $\tilde{f}$ and dividing
yields
\begin{equation}
\frac{\chi\left(t\right)}{\tilde{\chi}\left(t\right)}=\exp\left(\sum_{n=1}^{\infty}\frac{\kappa_{n}-\tilde{\kappa}_{n}}{n!}\left(it\right)^{n}\right),
\end{equation}
which, after rearrangement and taking the inverse Fourier transform,
yields
\begin{equation}
f\left(x\right)=\exp\left(\sum_{n=1}^{\infty}\frac{\kappa_{n}-\tilde{\kappa}_{n}}{n!}\left(-\frac{d}{dx}\right)^{n}\right)\tilde{f}\left(x\right).\label{eq:moment-expand}
\end{equation}
As with the Gram--Charlier series, the Edgeworth series is usually
presented for the Gaussian case $\tilde{f}=\Phi$. Although these
series are formally identical, they yield different partial sums when
truncated to a finite number of terms and hence have different convergence
properties. The Edgeworth form is generally considered more compact
than the Gram--Charlier series, as only the former is a true asymptotic
series.\cite{Blinnikov1998a,Stuart1994}

\subsection{Deriving the Gram--Charlier series from the Edgeworth series}

Rederiving the Gram--Charlier form from the Edgeworth series reveals
additional interesting relationships. One such relation follows from
the identity
\begin{equation}
\exp\left(\sum_{n=1}^{\infty}\frac{a_{n}}{n!}t^{n}\right)=\sum_{n=0}^{\infty}\frac{B_{n}(\left\{ a_{k}\right\} _{k=1}^{n})}{n!}t^{n},
\end{equation}
where $B_{n}$ is the complete Bell polynomial of order $n$ with
parameters $a_{1},\dots,a_{n}$.\cite{Bell1927} Setting $t=-d/dx$
gives immediately the differential operator
\begin{equation}
T\left(x\right)=\exp\left(\sum_{n=1}^{\infty}\frac{\left(\kappa_{n}-\tilde{\kappa}_{n}\right)}{n!}\left(-\frac{d}{dx}\right)^{n}\right)=\sum_{n=0}^{\infty}\frac{B_{n}\left(\left\{ \kappa_{k}-\tilde{\kappa}_{k}\right\} _{k=1}^{n}\right)}{n!}\left(-\frac{d}{dx}\right)^{n}.\label{eq:direct-series}
\end{equation}
We will call the last series of (\ref{eq:direct-series}) the direct
series of $T$.

Further specializing again to the Gaussian reference, we can use Rodrigues's
formula
\begin{equation}
He_{n}\left(x\right)\Phi\left(x\right)=\left(-1\right)^{n}\left(\frac{d}{dx}\right)^{n}\Phi\left(x\right),
\end{equation}
so that
\begin{equation}
T\left(x\right)\Phi\left(x\right)=\sum_{n=0}^{\infty}\frac{B_{n}\left(\left\{ \kappa_{k}-\tilde{\kappa}_{k}\right\} _{k=1}^{n}\right)}{n!}He_{n}\left(x\right)\Phi\left(x\right).
\end{equation}

The first few coefficients $c_{n}=B_{n}\left(\left\{ \kappa_{k}-\tilde{\kappa}_{k}\right\} _{k=1}^{n}\right)$
have been tabulated explicitly,\cite{Stuart1994} but to our knowledge
the relationship to the Bell polynomials have not been previously
discussed in the literature.

\subsection{Quantifying the effect of differing moments}

The Edgeworth series yields a useful result for error quantification.
If the first $k-1$ moments of two p.d.f.s $f$ and $\tilde{f}$ are
the same, but the $k$th moments differ, then the leading term in
the Edgeworth series is\begin{subequations}
\begin{align}
f\left(x\right) & =\tilde{f}\left(x\right)+\frac{\left(-1\right)^{k}B_{k}\left(\left\{ \kappa_{l}-\tilde{\kappa}_{l}\right\} _{l=1}^{k}\right)}{k!}\tilde{f}^{\left(k\right)}\left(x\right)+\mathcal{O}\left(\tilde{f}^{\left(k+1\right)}\right)\\
 & =\tilde{f}\left(x\right)+\frac{\left(-1\right)^{k}\left(\mu_{k}-\tilde{\mu}_{k}\right)}{k!}\tilde{f}^{\left(k\right)}\left(x\right)+\mathcal{O}\left(\tilde{f}^{\left(k+1\right)}\right).
\end{align}
\end{subequations}The second equality follows from the definition
of cumulants: the $k$th cumulant is a function of only the first
$k$ moments and can be written as $\kappa_{k}=\kappa_{k}\left(\mu_{1},\dots,\mu_{k}\right)=\mu_{k}+\dots$.

\subsection{The locus of differences in moments}

The first-order term in the preceding expansion is a quantitative,
asymptotic estimate of the difference between $f_{A+B}$ and $f_{A\star B}$.
The word expansion (\ref{eq:moment-expanded}) allows us to refine
the error analysis in terms of specific joint moments that contribute
to $\mu_{k}\left(A+B\right)-\mu_{k}\left(A\star B\right)$. Further
insight may be gained from the lattice sum approach pioneered by Wigner\cite{Wigner1955}
to interpret each term in (\ref{eq:moment-expanded}), being a trace
of a product of $k$ matrices, as a closed path with up to $k$ hops
as allowed by the structure of the matrices being multiplied.
\begin{example}
\label{example:tridiag-hopping}Consider $A$ and $B$ as in Example~\ref{example:tridiagonal},
which are partially free of degree 8 and whose discrepancy in the
eighth moments relative to complete freeness is solely in the term
$\left\langle \left(AB\right)^{4}\right\rangle $. $B$ is the adjacency
matrix of the one--dimensional chain $\cdot-\cdot-\cdots-\cdot-\cdot$
with $N$ nodes and periodic boundary conditions, and we can interpret
$\left\langle \left(AB\right)^{4}\right\rangle $ as the expected
sum of weights of particular paths on this lattice. These paths must
have exactly four hops, as $A$, being diagonal, does not permit hops,
whereas $B$, having nonzero entries only on the super- and sub-diagonals,
require exactly one hop either to the immediate left or the immediate
right. This gives rise to four different paths as illustrated in Figure~\ref{fig:tridiag-hopping}.
\begin{figure}[h]
\caption{\label{fig:tridiag-hopping} Paths contributing to the term $\left\langle \left(AB\right)^{4}\right\rangle $
in Example~\ref{example:tridiag-hopping}.}

\includegraphics{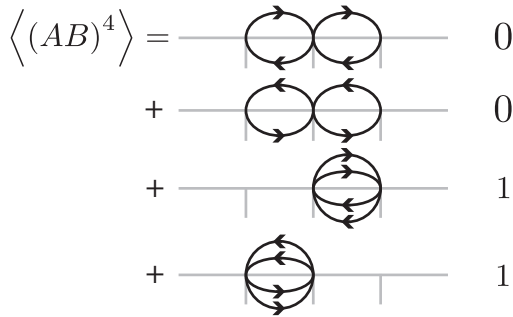}
\end{figure}
 We can show this by writing out the explicit matrix multiplication.
Writing the diagonal elements of $A$ as the i.i.d.\ standard Gaussians
$A_{ii}=g_{i}$ and using the Einstein implicit summation convention,
\begin{subequations}

\begin{align}
\left\langle \left(AB\right)^{4}\right\rangle = & \frac{1}{N}\mathbb{E}A_{i_{1}i_{2}}B_{i_{2}i_{3}}A_{i_{3}i_{4}}B_{i_{4}i_{5}}A_{i_{5}i_{6}}B_{i_{6}i_{7}}A_{i_{7}i_{8}}B_{i_{8}i_{1}}\\
= & \frac{1}{N}\mathbb{E}g_{i_{1}}g_{i_{3}}g_{i_{5}}g_{i_{7}}\delta_{i_{1}i_{2}}\left(\delta_{i_{2}-1,i_{3}}+\delta_{i_{2}+1,i_{3}}\right)\delta_{i_{3}i_{4}}\left(\delta_{i_{4}-1,i_{5}}+\delta_{i_{4}+1,i_{5}}\right)\nonumber \\
 & \qquad\times\delta_{i_{5}i_{6}}\left(\delta_{i_{6}-1,i_{7}}+\delta_{i_{6}+1,i_{7}}\right)\delta_{i_{7}i_{8}}\left(\delta_{i_{8}-1,i_{1}}+\delta_{i_{8}+1,i_{1}}\right)\\
= & \frac{1}{N}\mathbb{E}\left(g_{i}g_{i-1}g_{i}g_{i+1}+g_{i}g_{i+1}g_{i}g_{i-1}+g_{i}g_{i-1}g_{i}g_{i-1}+g_{i}g_{i+1}g_{i}g_{i+1}\right)\\
= & 2\mathbb{E}\left(g_{1}\right)^{2}\mathbb{E}\left(g_{1}^{2}\right)+2\mathbb{E}\left(g_{1}^{2}\right)^{2}=2.
\end{align}
\end{subequations}This simple example illustrates several important
concepts. First, if the underlying matrices can be interpreted as
adjacency matrices of graphs, these graphs are topologically significant
for the calculation of joint moments by controlling the number of
allowed returning paths in the lattice sum. The third equality makes
explicit use of this fact, and the numerical factor of 2 can be seen
as encoding the average degree of the underlying graph, which in this
case is the infinite one--dimensional chain. Second, the analyses
of joint moments of random matrices can be reduced to studying moments
and correlations of scalar matrix elements. While in this example
we assumed that the $g_{i}$s were uncorrelated, the calculation up
to the penultimate line is still valid in the general case where they
are correlated, and results like Wick's theorem\cite{Wick1950} or
Isserles's theorem\cite{Isserlis1916,Isserlis1918} can be applied
instead to finish the calculation. In an even more general setting,
other random variates other than standard Gaussians can be analyzed
with little difficulty so long as the required moments and correlations
exist. Third, these calculations can be performed for both finite
and infinite random matrices; for the former, there is even some information
about boundary conditions. If $B$ were replaced by $B^{\prime}$,
the implicit sum after the third equality excludes the some terms
for the boundary cases $i=1$ and $i=N$, since the paths would not
be able to travel beyond the edges of the lattice. As a result, the
numerical factors of 2 in the last line would be replaced by $2-2/N$
instead. In general, we expect boundary conditions to give rise to
corrections from the bulk behavior of order $\mathcal{O}\left(1/N\right)$.

To summarize this section, the notion of partial freeness unites two
disparate ideas in probability theory. First, the violation of free
independence in specific joint moments leads to asymptotic correction
factors to the density of states that show up as leading--order terms
in Edgeworth series expansions of the free convolution. These correction
factors have magnitudes that decay strongly with the lengths of the
words in question. Second, the coefficient of the correction also
encodes information about lattice sums over closed paths on random
graphs, whose topologies are encoded by the random matrices in question,
and can be related to quantities such as the average degree of a node
in the graph. Our generalization of the Edgeworth series to correct
for deviations from freeness (rather than to correct for nonnormality
in its classical usage) thus elucidates new connections between the
combinatorics of joint moments, sums over lattice paths on random
graphs, and asymptotic moment expansions of probability distributions.
\end{example}

\section{Computational implementation}

The relationship between joint moments and corrections to the density
of states is a particular feature of partial freeness which lends
itself naturally to numerical investigation. In this section, we sketch
how the characterization of partial freeness can be calculated in
an entirely automated fashion, by combining algorithms for enumerating
all joint moments of a given order with new algorithms using the results
above. Perhaps interestingly, partial freeness generates useful statistics
even when the analytic forms of the random matrices $A$ and $B$
are not known \emph{a priori}. We will treat this as a separate case
below.

\subsection{Symbolic computation of moments and joint moments}

If the analytic forms of the random matrices $A$ and $B$ are known
and can be multiplied analytically, a computer algebra system can
be used to calculate the necessary moments symbolically. The algorithm
for characterizing partial freeness then proceeds as follows:
\begin{enumerate}
\item Calculate the moments of $A$ and $B$ as well as the moments $\mu_{k}\left(A\boxplus B\right)$.
The $k$th moment of $A\boxplus B$ can be calculated by generating
all terms in the word expansion of $\left\langle \left(A+B\right)^{k}\right\rangle $
using Sawada's algorithm to generate all $\left(k,2\right)$--bracelets\cite{Sawada2001}.
\item For each word, check whether the relation (\ref{eq:free-defn}) required
by free independence holds. The first order $p$ for which this fails
is the degree of partial freeness.
\item Calculate the density of states $f_{A\boxplus B}$ using the $R$--transform
(\ref{eq:r-transform}) and its $p$th derivative $f_{A\boxplus B}^{\left(p\right)}$.
\item The leading--order correction to $f_{A\boxplus B}$ due to lack of
freeness is then $f_{A\boxplus B}+f_{A\boxplus B}^{\left(p\right)}\left(\mu_{p}-\tilde{\mu}_{p}\right)/p!$.
\end{enumerate}
To illustrate Step 2, we provide some \emph{Mathematica} code for
calculating the necessary moments and joint moments in Algorithm~\ref{alg:ma-joint-moments}.
The code also provides a function for calculating n.e.t.s of an arbitrary
joint moment or centered joint moment in terms of the distribution
of matrix elements. For simplicity, only the i.i.d.\ case of one
scalar probability distribution with moments $\left\{ m_{k}\right\} $
is illustrated, although this approach can be extended to more complicated
situations as necessary.

\begin{algorithm}[h]
\caption{\label{alg:ma-joint-moments}Mathematica code for calculating normalized
expected traces of joint matrix products and centered joint matrix
products for finite dimensional random matrices.}

\begin{lstlisting}[language=Mathematica]
NN = 100; (* Size of matrix *)

(* The following generates the map which
formally evaluates the expectation of the
G random variables assuming that they
are i.i.d. with vanishing mean. *)
MomentsOfG := Flatten[{
 Table[Subscript[G, j]^i -> Subscript[m, i],
  {i, 2, NN}, {j, 1, NN}],
 Table[Subscript[G, j] -> 0, {j, 1, NN}]
}];

ExpectationOfG[x_] := x /. MomentsOfG;

(* Normalized expected trace *)
AngleBracket[x_] := ExpectationOfG[
 Tr[x]/NN // Expand ];

(* centering operator *)
c[x_] := (x - AngleBracket[x] IdentityMatrix[NN])

(* Example 19 *)
A = DiagonalMatrix[Array[Subscript[G, #] &, NN]];
B = SparseArray[{Band[{1, 2}] -> 1}, {NN, NN}];
B[[1]][[-1]] = 1; (* Add circulant boundary *)
B = B + Transpose[B];

AngleBracket[c[A.A].c[B.B]]
(* Output: 0 *)

AngleBracket[MatrixPower[c[A].c[B], 4]]
(* Output: 2 Subscript[m, 2]^2 *)

(* Specialized to standard Gaussian Gs *)
GaussianG = Array[Subscript[m, #] ->
  Moment[NormalDistribution[], #] &, NN];
AngleBracket[MatrixPower[A+B], 4]] /. GaussianG

(* Output: 1099 *)
\end{lstlisting}
\end{algorithm}

\subsection{Numerical calculations on empirical samples}

The partial freeness formalism can also be used when the underlying
distributions of the random matrices $A$ and $B$ are unknown, but
when samples of each are available, e.g.\ from Monte Carlo simulations
or from empirical data. The algorithm for characterizing the partial
freeness of numerical samples is as follows:
\begin{enumerate}
\item Generate $t$ pairs of samples $\left\{ \left(A_{i},B_{i}\right)\right\} _{i=1}^{t}$
of $N\times N$ diagonalizable random matrices $A$ and $B$.
\item For each pair:

\begin{enumerate}
\item Calculate the exact eigenvalues of $A_{i}$ and $B_{i}$.
\item Calculate the eigenvalues of the sample $A_{i}+B_{i}$ of the exact
sum $A+B$. (This is for comparison purposes only and can be omitted.)
\item Calculate $n$ samples of the free convolution $f_{A\boxplus B}$
using the eigenvalues of $M_{i}=A_{i}+Q_{i}B_{i}Q_{i}^{\dagger}$
using a numerically generated Haar orthogonal (or unitary) matrix
$Q_{i}$.
\item Calculate $n$ samples from the classical convolution $f_{A\star B}$
using the eigenvalues of $L_{i}=A_{i}+\Pi_{i}^{-1}B_{i}\Pi_{i}$ using
a random permutation matrix $\Pi_{i}$.
\end{enumerate}
\item Calculate the first $2K$ moments of $A$ and $B$ as well as the
first $K$ moments of $\mu_{k}\left(A\boxplus B\right)$, $\mu_{k}\left(A\star B\right)$,
and $\mu_{k}\left(A+B\right)$.
\item Calculate the degree $k$ for which $A$ and $B$ are partially free
by testing for the smallest $k$ such that the moments of the free
convolution differ from the exact result, i.e. test the hypothesis
\[
\mu_{k}\left(A+B\right)\ne\mu_{k}\left(A\boxplus B\right).
\]

\item Using Sawada's algorithm,\cite{Sawada2001} enumerate all unique terms
$T_{j}=\left\langle A^{m_{1j}}B^{n_{1j}}\cdots A^{m_{k_{j}j}}B^{n_{k_{j}j}}\right\rangle $
in $\left\langle \left(A+B\right)^{k}\right\rangle $. For each term
$T_{j}$:

\begin{enumerate}
\item Calculate 
\[
T_{j}^{\left(cl\right)}=\left\langle A^{m_{1j}+\cdots+m_{k_{j}j}}\right\rangle \left\langle B^{n_{1j}+\cdots+n_{k_{j}j}}\right\rangle ,
\]
which would be its value expected from classical independence. Test
the hypothesis of equality $T_{j}=T_{j}^{\left(cl\right)}$.
\item Calculate the normalized expected trace of the centered term
\[
T_{j}^{\left(c\right)}=\left\langle \left(A^{m_{1j}}-\left\langle A^{m_{1j}}\right\rangle \right)\left(B^{n_{1j}}-\left\langle B^{n_{1j}}\right\rangle \right)\cdots\right\rangle ,
\]
which would be expected to vanish if $A$ and $B$ were truly free.
Test the hypothesis of equality $T_{j}^{\left(c\right)}=0$.
\end{enumerate}
\item Calculate the $k$th derivative $f_{A\boxplus B}^{\left(k\right)}$
using numerical finite difference.
\item Plot $f_{A+B}$, $f_{A\boxplus B}$ and $f_{A\boxplus B}+\left(\mu_{k}-\tilde{\mu}_{k}\right)/k!\cdot f_{A\boxplus B}^{\left(k\right)}$.
\end{enumerate}
This algorithm tests for partial freeness of degree $k\le K$, attempts
to identify the locus of discrepancy by testing all possible $\left(k,2\right)$-words,
and calculates the leading order correction term to the density of
states. The calculation of the classical convolution and exact density
of states are purely for comparative purposes and can be omitted.
In practice, we also account for sampling error in the hypothesis
tests in Steps 4 and 5 by calculating the standard error of each term
being tested, and evaluating the $p$--value for each such hypothesis.
For example, the standard error of the $k$th moment is
\begin{equation}
SE\left(\mu_{k}\right)=\sqrt{\frac{\mu_{2k}-\mu_{k}^{2}}{t}},
\end{equation}
and the standard error of a term $T_{j}=\left\langle \left(A^{m_{1j}}B^{n_{1j}}\cdots A^{m_{k_{j}j}}B^{n_{k_{j}j}}\right)\right\rangle $
in the expansion of $\left\langle \left(A+B\right)^{k}\right\rangle $
is
\begin{equation}
SE\left(T_{j}\right)=\sqrt{\frac{\left\langle \left(A^{m_{1j}}B^{n_{1j}}\cdots A^{m_{k_{j}j}}B^{n_{k_{j}j}}\right)^{2}\right\rangle -T_{j}^{2}}{t}}.
\end{equation}

The calculation of necessary standard errors require information up
to the $2K$th moment for calculating variances stemming from the
$K$th moment, which is why $2K$ moments of $A$ and $B$ are calculated
in Step 3. Alternatively, other measures of statistical fluctuation
could be used, such as bootstrap or jackknife errors.

The \emph{Supplementary Information} includes an implementation of
the algorithm described in this section for analyzing empirical samples
of random matrices that is written in MATLAB. In numerical tests of
this algorithm, we observe the expected $O\left(1/\sqrt{Nt}\right)$
rate of convergence in the word values with the number of eigenvalues
$Nt$. Thus in practical numerical studies where $N$ and $t$ can
be controlled, we recommend that for maximum numerical efficiency
that $N$ be set only as large as necessary to minimize finite-size
effects, and $t$ be taken as large as necessary to ensure numerical
convergence, as the diagonalization of typical matrices is superlinear
in $N$.

\section{Summary}

Partial freeness is a relationship between random matrices that brings
together ideas from various aspects of probability theory. First,
the notion of asymptotic freeness arises naturally as a special limiting
case of partial freeness for infinite--dimensional matrices, but unlike
the former, partial freeness is still well--defined for arbitrary
diagonal random or deterministic matrices of finite or infinite dimensions.
Second, partial freeness allows for deviations from asymptotic freeness
to be quantified in terms of well--defined asymptotic corrections
to quantities such as the empirical density of states. These asymptotic
corrections generalize the notions of Gram--Charlier and Edgeworth
series which arise from classical probability in the study of deviations
from non--Gaussianity. Third, the organization of joint moments by
words of a given length reveals new combinatorial structure, which
to our knowledge, has not been elucidated in the context of free probability
before. The enumeration of joint moments evokes the combinatorics
of necklaces, which also shows how the ideas in this paper generalize
straightforwardly to multiple additive free convolutions: the $k$
parameter of the necklaces in Definition~\ref{def:necklace} counts
the number of matrices whose sum $M=A+B+\dots$ is being investigated.

We have also demonstrated that partial freeness is a theoretically
interesting abstract relation between random matrices, but it also
comes with a statistical framework which can be tested in numerical
computations in a practical manner. Partial freeness organizes clearly
the relationships between the joint moments of random matrices and
the moments and correlations of the scalar random variables in their
matrix elements. Additionally, partial freeness can be tested for
statistically using purely empirical data, without resorting to any
model for the random matrices in question. These ideas can be stated
in algorithmic form and thus we expect partial freeness to be useful
both theoretically and in practical numerical applications. We are
currently exploring how the theoretical ideas brought together by
partial freeness can be used to construct new computational statistical
tools.

\thanks{We acknowledge funding from NSF SOLAR Grant No.~1035400. A.E.\ acknowledges
additional funding from NSF DMS Grant No.~1016125. We gratefully
acknowledge useful discussions with D. Shlyakhtenko (UCLA), N. Raj
Rao (Michigan) and A. Su\'arez (Univ. Aut\'onoma Madrid) that have
led us to pursue this avenue of investigation. We thank M. Welborn
(MIT) and E. Hontz (MIT) for graphics assistance with Figure~\ref{fig:tridiag-hopping}
and Figure~\ref{fig:tridiag-hopping} respectively.}

\bibliographystyle{focm}
\bibliography{partial-freeness}

\end{document}